\documentclass[12pt,reqno]{article}
\usepackage{amssymb,amscd,amsmath,amsthm,color}
\newtheorem{theorem}{Theorem}[section]
\newtheorem{lemma}[theorem]{Lemma}
\newtheorem{cor}[theorem]{Corollary}

\usepackage{enumerate,amssymb}
\theoremstyle{definition}

\newtheorem{example}[theorem]{Example}

\theoremstyle{remark}

\numberwithin{equation}{section}

\def\bM{\mathbb{M}}

\def\bR{\mathbb{R}}

\def\bE{\mathbb{E}}
\def\bB{\mathbb{B}}

\begin{document}
\baselineskip=15pt

\title{Positive linear maps on Hilbert space operators and noncommutative $L_p$ spaces}

\author{Jean-Christophe Bourin\footnote{Funded by the ANR Projet (No.\ ANR-19-CE40-0002) and by the French Investissements
 d'Avenir program, project ISITE-BFC (contract ANR-15-IDEX-03).
}\,  and Jingjing Shao\footnote{
The research is partially supported by the National Natural Science Foundation of
China No. 11701255}}

\date{ }

\maketitle

\vskip 10pt\noindent
{\small
{\bf Abstract.} We extend some inequalities for normal matrices and positive linear maps related to the Russo-Dye theorem.
The results cover the case of some positive linear maps $\Phi$ on a von Neumann algebra ${\mathcal{M}}$ such that $\Phi(X)$ is unbounded for all nonzero $X\in{\mathcal{M}}$.

}

\section{Some matrix inequalities}

We aim to study positive linear maps taking values in some spaces of not necessarily bounded Hilbert space operators.
This study is driven by some recent matrix inequalities   established in \cite{BLpos} and \cite{BL-russodye}. 
The space of $n\times n$ complex matrices is denoted by $\bM_n$ and its positive semi-definite cone by $\bM_n^+$.
The identity, in any algebra through the text, is denoted by $I$.
A linear map $\Phi: \mathrm{M}_n\rightarrow \mathrm{M}_m$ is  positive
if $\Phi(\mathrm{M}_n^+)\subset \mathrm{M}_m^+$.

Let $\Phi: \mathrm{M}_n\rightarrow \mathrm{M}_m$ be a positive linear map and let $N\in\bM_n$ be normal. Then
there exists a unitary $V\in\mathrm{M}_m$ such that
\begin{equation*}\label{comp1}
|\Phi(N)|\leq \frac{\Phi(|N|)+V\Phi(|N|)V^*}{2}.
\end{equation*}
and
\begin{equation*}\label{comp2}
|\Phi(N)|\leq \Phi(|N|)+\frac{1}{4} V\Phi(|N|)V^*.
\end{equation*}
These two inequalities and
several consequences are proved in \cite{BLpos}, \cite{BL-russodye}.  As an application for the Schur product of two normal matrices $A, B\in\mathrm{M}_n$, one may infer that
\begin{equation*}\label{inequality for Schur product}
|A\circ B|\leq|A|\circ|B|+\frac{1}{4}V(|A|\circ|B|)V^*
\end{equation*}
for some unitary $V\in\mathrm{M}_n$, where the constant $1/4$ is optimal.
Another interesting consequence is the following improvement of the Russo-Dye theorem stating that every positive linear map attains its norm at the identity: if $Z\in\bM_n$ is a contaction, then
\begin{equation*}\label{RDimprov}
 |\Phi(Z)|  \le \frac{\Phi(I) +V\Phi(I)V^*}{2}.
\end{equation*}
Applying this to the Schur product with $S\in\bM_n^+$ yields some exotic eigenvalue inequalities such as
\begin{equation*}\label{eqschur}
\lambda_3(|S\circ Z|) \le \delta_{2}(S)
\end{equation*}
where $\lambda_3(\cdot)$ stands for the third largest eigenvalue, and $\delta_2(\cdot)$ for the second largest diagonal entry.

In the next section we shall extend these inequalities to the setting of Hilbert space operators, with a special emphasis on notions that do not exist on the matrix case, such as
hyponormal or semi-hyponormal operators, and the Calkin theory of operator ideals.

Sections 3 and 4 are devoted to  unbounded operators and unbounded positive linear maps. The correct setup is that of
$\tau$-measurable operators affiliated to semi-finite von Neumann algebras and positive linear maps, continuous with respect to the measure topologies. We recall these notions in Section 3, with some natural examples of positive linear maps,  some of  which can be regarded are purely unbounded. The proofs for measurable operators are rather different than those for matrices.
In particular, contrarily to the matrix case, we will not use the geometric mean.

\section{Positive maps taking values in $\bB(\mathcal{H})$}

Denote by $\mathcal{A}$ a unital  $C^*$-algebra acting on an infinite dimensional  separable Hilbert space
$\mathcal{H}$, and let  $\bB(\mathcal{H})$   stand for the set of all bounded linear operators on $\mathcal{H}$.
Let $A,B\in\bB(\mathcal{H})$ be positive and invertible. Their geometric mean is defined as
$$
A\#B :=A^{1/2}(A^{-1/2}BA^{-1/2})^{1/2}A^{1/2}
$$
As in the matrix case, the geometric cannot be extended by continuity (in norm, or even in the strong operator topology) to positive, noninvertible operators, however, the natural definition for
 positive, not necessarily invertible operators, is the strong limit
$$
A\#B :={\mathrm{sot }}\,\lim_{r \to 0^+} (A+rI)\#(B+rI).
$$
We then have the arithmetic-geometric mean inequality
$$
A\# B \le \frac{A+B}{2}
$$
and, replacing $A$ by $2A$, $B$ by $(1/2)B$,
$$
A\# B \le A +\frac{1}{4}B.
$$
Hence, the next theorem contains several arithmetic means inequalities.
\vskip 5pt
\begin{theorem}\label{thmBH}     Let $\Phi: {\mathcal{A}}\to \bB(\mathcal{H})$   be  a positive linear map and let $N\in {\mathcal{A}}$  be normal. Then, there exists  a partial isometry $V\in\bB(\mathcal{H})$ such that
$$
   |\Phi(N)|  \le \Phi(|N|) \#V\Phi(|N|)V^*.
$$
\end{theorem}

\vskip 5pt
\begin{proof} The proof is  the same as in the matrix case, see the proof of \cite[Eq. (2.1)]{BL-russodye}.
\end{proof}

\vskip 5pt
\begin{cor} \label{corideal1}  Let $\Phi: {\mathcal{A}}\to \bB(\mathcal{H})$   be  a positive linear map, and let $N\in{\mathcal{A}}$ be normal. If $\Phi(|N|)$ belongs to an ideal ${\mathcal{I}}\subset\bB(\mathcal{H})$, then so does $\Phi(N)$.
\end{cor}

\vskip 5pt
\begin{proof} From the arithmetic-geometric mean inequality
$$
   |\Phi(N)|  \le \Phi(|N|) \#V\Phi(|N|)V^* \le \frac{ \Phi(|N|) +V\Phi(|N|)V^*}{2}
$$
we infer that $|\Phi(N)|$, and so $\Phi(N)$, belong to ${\mathcal{I}}$.
\end{proof}

\vskip 5pt
To delete the normality asssumption on $N$ in Corollary \ref{corideal1}, we first give one more consequence of Theorem \ref{thmBH}.

\vskip 5pt
\begin{cor}\label{corcons2}     Let $\Phi: {\mathcal{A}}\to \bB(\mathcal{H})$   be  a positive linear map and let $X\in {\mathcal{A}}$. Then, there exists  a partial isometry $V\in\bB(\mathcal{H})$ such that
$$
   |\Phi(X\pm X^*)| \le \frac{\Phi(|X|+|X^*|) + V\Phi(|X|+|X^*|)V^*}{2}.
$$
\end{cor}

\vskip 5pt
\begin{proof} Let $\Psi:\bM_2({\mathcal{A}})\to \bB(\mathcal{H})$ be defined as
$$
\Psi\left( \begin{pmatrix} A&B \\ C&D\end{pmatrix}\right)=\Phi(A+B+C+D).
$$
Since
$$
A+B+C+D=\begin{pmatrix} I&I\end{pmatrix}\begin{pmatrix} A&B \\ C&D\end{pmatrix}\begin{pmatrix} I\\ I\end{pmatrix},
$$
 $\Psi$ is a positive map. Applying Theorem  \ref{thmfinite}   to this map with the normal (Hermitian) operator in $\bM_2({\mathcal{A}})$
$$
\begin{pmatrix} 0&X \\ X^*&0\end{pmatrix}
$$
yields the result with the $+$ sign. Replacing $X$ by $iX$ yields the $-$ sign case.
\end{proof}

\vskip 5pt
We are now in a position to delete the normality assumption in Corollary \ref{corideal1}.

\vskip 5pt
\begin{cor} \label{corideal2}  Let $\Phi: {\mathcal{A}}\to \bB(\mathcal{H})$   be  a positive linear map, and let $A\in{\mathcal{A}}$ be invertible. If $\Phi(|A|)$ belongs to an ideal ${\mathcal{I}}\subset\bB(\mathcal{H})$, then so does $\Phi(A)$.
\end{cor}

\vskip 5pt
Before giving the proof, we show that the invertibility assumption on $A$ is crucial, even for completely positive linear maps. Let ${\mathcal{A}}=\bB(\mathcal{H})$  and pick the infinite direct sum copies
$$
A:=\bigoplus^{\infty}
\begin{pmatrix}
0&0 \\
1&0
\end{pmatrix}
$$
so that
$$
|A|=\bigoplus^{\infty}
\begin{pmatrix}
1&0 \\
0&0
\end{pmatrix}.
$$
Consider the Schur multiplier $\Phi(X)=S\circ X$ with
$$
S=\bigoplus_{n=1}^{\infty}
\begin{pmatrix}
1/n^2&1/n \\
1/n& 1
\end{pmatrix}.
$$
Hence,
$$
S\circ |A| = \bigoplus_{n=1}^{\infty}
\begin{pmatrix}
1/n^2&0 \\
0& 0
\end{pmatrix}
$$
belongs to the Trace-class ideal, while
$$
S\circ A = \bigoplus_{n=1}^{\infty}
\begin{pmatrix}
0&0 \\
1/n& 0
\end{pmatrix}
$$
does not belong to the Trace-class ideal.

We turn to the proof of Corollary \ref{corideal2}.

\vskip 5pt
\begin{proof} Since $A$ is invertible, so is $A^*$. Let
$$
\| A^*A^{-1}\|_{\infty} =\nu.
$$
Then,
$$
A^{-1*}AA^*A^{-1} \le \nu^2 I,
$$
equivalently,
$$
AA^* \le \nu^2A^*A,
$$
and so,
\begin{equation}\label{hypo}
|A^*| \le \nu|A|.
\end{equation}
Now, observe that
$$
\Phi(A)=\frac{\Phi(A+A^*)}{2} + \frac{\Phi(A-A)^*}{2}
$$
Now, from Corollary \ref{corcons2} and \eqref{hypo} we  have a partial isometry $V_1\in\bB(\mathcal{H})$ such that
$$
\left| \Phi(A+A^*)\right| \le \frac{\Phi(|A|+|A^*|) +V_1 \Phi(|A|+|A^*|)V_1^*}{2}\le (1+\nu)\frac{\Phi(|A|) +V_1 \Phi(|A|)V_1^*}{2}.
$$
Hence,  $\Phi(A+A^*) \in\mathcal{I}$. Similarly  $\Phi(A-A^*) \in\mathcal{I}$, and so
 $\Phi(A) \in\mathcal{I}$.
\end{proof}

\vskip 5pt
We wish to extend Theorem \ref{thmBH} for normal operators to a larger class of operators. Recall  the following classical extensions of normal operators:
$$
{\text{Normal }} \ \subset \   {\text{Quasinormal }}\  \subset \ {\text{Subnormal }}\  \subset \ {\text{Hyponormal}}\  \subset \ {\text{Semi-hyponormal }}
$$
where the inclusions are strict and the larger class of semi-hyponormal operators $S$ is defined by the condition $|S^*|\le |S|$. To this end we first need the following consequence of  Theorem \ref{thmBH}.

\vskip 5pt
\begin{cor}\label{corRD1}     Let $\Phi: {\mathcal{A}}\to \bB(\mathcal{H})$   be  a positive linear map and let $Z\in {\mathcal{A}}$  be a contraction. Then, there exists  a partial isometry $V\in\bB(\mathcal{H})$ such that
$$
   |\Phi(Z)|  \le \Phi(I) \#V\Phi(I)V^*.
$$
\end{cor}

\vskip 5pt
This result is a far extension of the famous Russo-Dye Theorem asserting that every positive linear map $\Phi $ on a unital $C^*$-algebra attains its norm at the identity,
$$
\|\Phi(Z)\|_{\infty} \le \|\Phi(I)\|_{\infty}
$$
for all contractions $Z$. In the matrix setting, we refer to
\cite{BL-russodye} for more general results than Corollary \ref{corRD1} and several applications. The proof for operators is exactly the same and is given here for sake of convenience.

\vskip 5pt
\begin{proof} We may dilate $Z$ into a unitary $U\in\bM_2({\mathcal{A}})$, for instance with Halmos,
$$
U=\begin{pmatrix}
Z&-\sqrt{I-ZZ^*} \\
\sqrt{I-Z^*Z}& Z^*
 \end{pmatrix}
$$
Now, let $\Psi:\bM_2({\mathcal{A}})\to\bB({\mathcal{H}})$ be defined as
$$
\Psi\left( \begin{pmatrix} A&B \\ C&D\end{pmatrix}\right)=\Phi(A).
$$
Applying Theorem \ref{thmBH} to $\Psi$ and $U$, we have
$$
|\Phi(Z)|=|\Psi(U)|\le \Psi(|U|) \#V\Psi(|U|)V^*= \Phi(I) \#V\Phi(I)V^*
$$
for some unitary $V\in\bB({\mathcal{H}})$.
\end{proof}

\vskip 5pt
\begin{theorem}\label{thmsemi}     Let $\Phi: {\mathcal{A}}\to \bB(\mathcal{H})$   be  a positive linear map and let $S_1,S_2,\ldots,S_m\in {\mathcal{A}}$  be semi-hyponormal. Then, there exists a partial isometry $V\in\bB(\mathcal{H})$ such that
$$
   \left|\Phi\left(\sum_{k=1}^m S_k\right)\right|  \le \Phi\left(\sum_{k=1}^m |S_k|\right)\#V\Phi\left(\sum_{k=1}^m |S_k|\right)V^*.
$$
\end{theorem}

\vskip 5pt
\begin{proof} By considering the semi-hyponormal operator
$$
S:=S_1\oplus\cdots\oplus S_m\in\oplus^m{\mathcal{A}}
$$
and the positive linear map $\Psi:\oplus^m{\mathcal{A}}\to\ \bB(\mathcal{H})$ defined as
$$
\Psi\left( X_1\oplus\cdots\oplus X_m\right) =\Phi\left(\sum_{i=1}^m X_i\right)
$$
it suffices to prove the theorem for $m=1$, i.e, for the single semi-hyponormal operator $S$ and the map $\Psi$.

Consider the map $\Lambda:\oplus^m{\mathcal{A}}\to\bB({\mathcal{H}})$,
$$
\Lambda(X) =\Psi\left(|S|^{1/2}X(S|^{1/2}\right).
$$
Observe that
\begin{equation}\label{eqt3}
\Lambda(I)=\Psi(|S|), \quad {\mathrm{and}} \quad \Lambda(Y)=\Psi(S)
\end{equation}
where
$$Y=|S|^{-1/2}S|S|^{-1/2}
$$
and $|S|^{-1}$ stands for the generalized inverse.
Thanks to the polar decomposition $S=|S^*|U=U|S|$, we  have
$$Y=|S|^{-1/2}|S^*|^{1/2}U|S|^{1/2}|S|^{-1/2}
$$
and the semi-hyponormality assumption on $S$ entails that $|S|^{-1/2}|S^*|^{1/2}$ and, of course, $|S|^{1/2}|S|^{-1/2}$, the support projection of $|S|$, are two contractions. Therefore $Y$ is a contraction too. Applying Corollary \ref{corRD1} to $Y$ and $\Lambda$ yields
\begin{equation*}
|\Lambda(Y)| \le \Lambda(I) \#V\Lambda(I)V^*
\end{equation*}
for some partial isometry $V\in\bB(\mathcal{H})$. Coming back to \eqref{eqt3} we get
$$
|\Psi(S)| \le \Psi(|S|)\#V\Psi(|S|)V^*
$$
which completes the proof.
\end{proof}

\vskip 5pt
We close this section by some application to Cartesian decomposition, following \cite[Corollary 3.5]{BL-russodye}.

\vskip 5pt
\begin{cor} \label{corgeo1}  Let $\Phi: {\mathcal{A}}\to \bB(\mathcal{H})$   be  a positive linear map and let $Z\in {\mathcal{A}}$  with Cartesian decomposition $Z=X+iY$. Then, there exists  a partial isometry $V\in\bB(\mathcal{H})$ such that
$$
   |\Phi(Z)|  \le \Phi(|X|+|Y|) \#V\Phi(|X|+|Y|)V^*.
$$
\end{cor}

\vskip 5pt
A special case of this corollary with the identity map $\Phi(Z)=Z$ combined with the arithmetic-geometric mean inequality reads as
\begin{equation*}
|Z|  \le \frac{|X|+|Y| +V(|X|+|Y|)V^*}{2}.
\end{equation*}
Corollary \ref{corgeo1} also  yields an inequality for the essential norm,
$$
   \|\Phi(Z)\|_{ess}  \le \|\Phi(|X|+|Y|)\|_{ess}.
$$

\section{Positive linear maps on $\tau$-measurable operators}

In Section 4, we will extend the results of Section 2 to unbounded operators. The correct framework consists in measurable operators affiliated to a semifinite von Neumann algebra $\mathcal{M}$, acting on a separable Hilbert space, with a
faithful normal semifinite trace $\tau$. In this short section we provide an example of a positive linear map which maps every nonzero positive operator in  $\mathcal{M}$ to an unbounded operator. Such a map cannot be norm continuous; for measurable operators, the notion of convergence in measure confers the good topology. We refer the reader to the Fack-Kosaki survey
\cite{FK} for a nice detailed survey on this theory.

We recall the notion of convergence in measure.
Let $\overline{\mathcal{M}}$ denote the set of $\tau$-measurable operators affiliated with $\mathcal{M}$ and $\overline{\mathcal{M}}^+$ the positive cone of $\overline{\mathcal{M}}$. The
spectral scale of $A\in\overline{\mathcal{M}}^{+}$ is defined as
\begin{equation}\label{sp-scale}
\lambda_t(A):=\inf\{s\in\bR:\tau({\mathbf{1}}_{(s,\infty)}(A))\le t\},
\qquad t\in(0,\tau(I)),
\end{equation}
where ${\mathbf{1}}_{(s,\infty)}(A)$ is the spectral projection of $A$ corresponding to
$(s,\infty)$. The generalized $s$-numbers of $X\in \overline{\mathcal{M}}$ is
$\mu_t(X):=\lambda_t(|X|)$, $t\in(0,\tau(I))$. A sequence $\{X_n\}$ in $\overline{\mathcal{M}}$ converges in measure to $X\in\overline{\mathcal{M}}$
if for all $t>0$, $\mu_t(X_n)\to \mu_t(X).$

Hence, in case of ${\mathcal{M}}=\bB({\mathcal{H}})$ we merely have $\overline{\mathcal{M}}={\mathcal{M}}$, and the convergence in measure coincides with the norm convergence. For a diffuse von Neumann algebra, the situation is more interesting, and it makes sense to consider unbounded positive linear maps which are continuous with respect to the measure topologies.

\vskip 5pt
\begin{example}\label{ex1}
Let $\{Z_i\}_{i=1}^m\subseteq \overline{\mathcal{M}}$ and define the  map $\Phi:\overline{\mathcal{M}}\to\overline{\mathcal{M}}$ by
$$\Phi(x)=\sum_{i=1}^mZ_i^*XZ_i.$$
Since $\overline{\mathcal{M}}$ is a complete topological (metrizable) algebra, $\Phi$ is the most basic and natural example of positive map, continuous with respect to the measure topology. If some of the weights $Z_i$ are not bounded, then this map is not norm continuous.
\end{example}

\vskip 5pt
The space $\overline{\mathcal{M}}$ is often denoted by $L_0({\mathcal{M}})$ to recall the continuous embeddings of the  noncommutative $L_p({\mathcal{M}})$ spaces, $0<p<\infty$, into
$\overline{\mathcal{M}}$. However, some natural positive linear maps cannot be defined
on the whole space $\overline{\mathcal{M}}$, for instance, the conditional expectation onto a subalgebra cannot be defined on $L_q({\mathcal{M}})$, $0\leq q<1$.

\vskip 5pt
\begin{example}\label{ex2} Let $\mathcal{M}\subset\mathcal{N}$ be two type ${\mathrm{II}}_1$ factors, and let $\bE:\mathcal{M}\to \mathcal{N}$ be the (trace preserving) conditional expectation from $\mathcal{M}$ to its subfactor $\mathcal{N}$. This map  is  continuous for both the norm topologies and the $L_1$-norm topologies, in fact $\bE$ is a contractive map from $L_p({\mathcal{M}})$ to $L_p({\mathcal{M}})$, $1\le p\le\infty$, as in the standard commutative case. Let $\{Z_i\}_{i=1}^m\subseteq \overline{\mathcal{M}}$ and define the map $\Phi: L_1({\mathcal{M}})\to\overline{\mathcal{M}}$ by
$$\Phi(x)=\sum_{i=1}^mZ_i^*\bE(X)Z_i.$$
The map $\Phi$ is continuous (for the $L_1$ norm and the measure topology).
\end{example}

\vskip 5pt
\begin{example}\label{ex3} Let $\mathcal{M}$ be a type ${\mathrm{II}}_1$ factor acting on the Hilbert space ${\mathcal{H}}$, let $\{ h_n\}_{n=1}^{\infty}$ be a sequence of vectors dense in the unit sphere of  ${\mathcal{H}}$, and define a positive linear map on $\Phi:{\mathcal{M}}\to \overline{\mathcal{M}}$ by
$$\Phi(X)=\sum_{n=1}^{\infty} \langle h_n, Xh_n\rangle Y_n$$
where $Y_n\in\overline{\mathcal{M}}^+\setminus {\mathcal{M}}^+$ satisfy $Y_nY_k=0$ for $n\neq k$.
Then  $\Phi$ is  continuous (for the  norm on  $\mathcal{M}$  and the measure topology on $\overline{\mathcal{M}}$) and $\Phi(X)\in\overline{\mathcal{M}}\setminus {\mathcal{M}}$ for all nonzero $X\in{\mathcal{M}}$.
\end{example}

\vskip 5pt
\begin{example}\label{ex4} Let $\mathcal{M}$ be a type ${\mathrm{II}}_1$ factor acting, let $\{ H_n\}_{n=1}^{\infty}$ be a sequence in the unit ball of ${\mathcal{M}}$ dense for the $*$-weak topology $\sigma({\mathcal{M}}, L_1({\mathcal{M}}))$
and define a positive linear map on $\Phi:L_1({\mathcal{M}})\to \overline{\mathcal{M}}$ by
$$\Phi(X)=\sum_{n=1}^{\infty} \tau(|H_n|X) Y_n$$
where $Y_n\in\overline{\mathcal{M}}^+\setminus L_1({\mathcal{M}})^+$ satisfy $Y_nY_k=0$ for $n\neq k$.
Then  $\Phi$ is  continuous (for the $L_1$ norm   and the measure topology) and $\Phi(X)\in\overline{\mathcal{M}}\setminus L_1({\mathcal{M}})$ for all nonzero $X\in L_1({\mathcal{M}})$.
\end{example}

\section{Maps taking values in a type ${\mathrm{II}}$ factor}

There is no suitable definition of the geometric mean for positive $\tau$-measurable operators in a diffuse semifinite von Neumann algebra, and though some candidates might be considered, their basic properties are still not understood. We thanks Fumio Hiai for interesting discussions \cite{Hiai} on the current state of art and for showing us his forthcoming work on this  topic. It is worth  mentioning that Hiai is able to define the geometric mean in the noncommutative $L^p$ spaces, and one may expect that this could be extended to the set of all measurable operators.

Another difficulty for unbounded positive linear maps defined on a commutative domain is the lack of a Stinespring's lemma ensuring that the map is completely positive.

Therefore the main proof in this section is rather different from that in the matrix or bounded operator case.

To extend the results of Section 2 to the setting of $\tau$-measurable operators, we cannot use the geometric mean but we are still able to establish interesting  inequalities involving means in unitary (or partial isometry) orbits.

Let $\mathcal{M}$ denote a semifinite von Neumann algebra and let ${\mathcal{F}}_1$ be a type ${\mathrm{II}}_1$ factor with a normalized trace, i.e., taking the value 1 on the identity.

The
{\it spectral scale} of a self-adjoint operator $A\in\overline{\mathcal{F}}_1$ is defined as
\begin{equation}\label{sp-scale}
\lambda_t(A):=\inf\{s\in\bR:\tau({\mathbf{1}}_{(s,\infty)}(A))\le t\},
\qquad t\in(0,1),
\end{equation}
where ${\mathbf{1}}_{(s,\infty)}(A)$ is the spectral projection of $A$ corresponding to
$(s,\infty)$. The function $t\mapsto\lambda_t(A)$
is non-increasing and right-continuous. We may write $\lambda_0(A)$ and
$\lambda_{1}(A)$ for $\lim_{t\searrow0}\lambda_t(A)$ and
$\lim_{t\nearrow 1}\lambda_t(A)$, respectively, (which are the maximal and minimal
spectra of $A$ when $A$ is bounded). Note that the generalized $s$-numbers of $X\in\overline{\mathcal{F}}$ is
$\mu_t(X)=\lambda_t(|X|)$, $t\in(0,1)$.

The following two lemmas belong to the folklore.

\vskip 5pt
\begin{lemma}\label{lemma-sp-dominance}
For two self-adjoint operators  $A,B\in\overline{\mathcal{F}}_1$, the spectral dominance relation
$\lambda_t(A)\le\lambda_t(B)$ for all $t\in(0,1)$ holds if and only if for every $\varepsilon>0$ there
exists a unitary $U\in{\mathcal{F}}$ such that $UAU^*\le B+\varepsilon I$.
\end{lemma}

\vskip 5pt
\begin{lemma}\label{lemma-ae-convergence}
Let $\{A_n\}\subset\overline{\mathcal{F}}_1$ be a sequence of self-adjoint operators converging in measure to $A$. Then,
for all points of continuity $t$ of $s\mapsto \lambda_s(A)$, hence almost everywhere,
$$
\lim_{n\to \infty}\lambda_t(A_n)= \lambda_t(A).
$$
\end{lemma}

\vskip 5pt For positive linear maps taking value into the space of measurable operators affiliated to a type ${\mathrm{II}}_1$ factor, the following theorem holds.

\vskip 5pt
\begin{theorem}\label{thmfinite}     Let $\Phi: L_p({\mathcal{M}})\to \overline{\mathcal{F}}_1$ be a continuous positive linear map, let $N\in L_p({\mathcal{M}})$ be normal, $\beta> 0$ and $\varepsilon>0$. Then, there exists a unitary $V\in{\mathcal{F}}_1$ such that
$$
   |\Phi(N)|  \le \beta\Phi(|N|) +\frac{1}{4\beta}V\Phi(|N|)V^* +\varepsilon I.
$$
\end{theorem}

\vskip 5pt
The continuity assumption refers to the natural topology on $L_p({\mathcal{M}})$ and $\overline{\mathcal{F}}_1$. The most important cases are
$p=0, 1, \infty$, i.e, when $\Phi$ is defined on $\overline{\mathcal{M}}$,
$L^1(M)$, or ${\mathcal{M}}$.

 \vskip 5pt
\begin{proof} We write the proof for the case $p=0$, the other cases being quite similar. We start with the case $N\in {\mathcal{M}}$ and  $\Phi(I)\in{\mathcal{F}}_1$. Since any $X\in{\mathcal{M}}$ is a linear combination of four positive elements, the assumption $\Phi(I)$ is bounded is equivalent to the fact that
$\Phi$ induces a bounded positive linear map from ${\mathcal{M}}$ into ${\mathcal{F}}$. By confining $\Phi$ to the abelian unital $C^*$-algebra spanned by $N$, we may  assume that $\Phi$ is completely positive thanks to Stinespring's lemma.
Since in the algebra $\bM_2({\mathcal{M}})$,
\begin{equation}\label{eqfund}
\begin{pmatrix}
|N|& N\\
N^*& |N|
\end{pmatrix} \ge 0.
\end{equation}
we then infer that, in the algebra $\bM_2({\mathcal{F}})$,
\begin{equation}\label{eqfund}
\begin{pmatrix}
\Phi(|N|)& \Phi(N)\\
\Phi(N^*)& \Phi(|N|)
\end{pmatrix} \ge 0.
\end{equation}
Now, let $V^*$ be the unitary part in the polar decomposition $\Phi(N)=V^*|\Phi(N)|$. We have, for all $\alpha>0$,
\begin{equation*}
\begin{pmatrix} \alpha^{-1/2} V & -\alpha^{1/2} I \end{pmatrix}
\begin{pmatrix}
\Phi(|N|)& \Phi(N)\\
\Phi(N^*)& \Phi(|N|)
\end{pmatrix}
\begin{pmatrix} \alpha^{-1/2} V^*\\ -\alpha^{1/2} I \end{pmatrix}
\ge 0.
\end{equation*}
Equivalently,
$$
   |\Phi(N)| \le \frac{\alpha\Phi(|N|) +\alpha^{-1}V\Phi(|N|)V^*}{2}
$$
and setting $\beta=2\alpha$ yields
the inequalities of the theorem with $\varepsilon=0$.

We turn to the the general case.
Let
$$
\Phi(I) =\int_{0}^{\infty} \lambda \,{\mathrm{d}}E(\lambda)
$$
and define
$$
K_n=\int_{0}^{\infty} g_n(\lambda) \,{\mathrm{d}}E(\lambda)
$$
with $g_n(\lambda)=1$ for $\lambda\le n$ and $g_n(\lambda)=\lambda^{-1}$ for $\lambda> n$. Thus $K_n\in{\mathcal{F}}^+$ with $\| K_n\|_{\infty}=1$. Then define the positive linear maps
$$
\Phi_n(T) :=K_n^{1/2}\Phi(T)K_n^{1/2}
$$
and observe that these maps are bounded on ${\mathcal{M}}$ as $\|\Phi_n(I)\|_{\infty}\le n$.
Since $\{K_n\}$ is a (bounded) sequence converging in measure (denoted by $\to$) to $I$, we infer that $K_n^{1/2}Y_nK_n^{1/2}\to Y$ for any sequence $\{Y_n\}$ in $\overline{\mathcal{M}}$ such that $Y_n\to Y$. Hence $\Phi_n(T_n)\to \Phi(T)$ for any sequence $\{T_n\}$ in $\overline{\mathcal{M}}$ such that $T_n\to T$.

Now, let $N\in\overline{\mathcal{M}}$ be normal and pick a normal sequence $\{N_n\}\subset {\mathcal{M}}$ such
that $N_n\to N$. By the first step of the proof, we have
$$
   |\Phi_n(N_n)|-\beta \Phi(|N_n|)\le \frac{1}{4\beta}V_n\Phi(|N|)V_n^*
$$
where $V_n^*$ is the unitary part in the polar decomposition $\Phi(N_n)=V_n^*|\Phi(N_n)|$. This entails that, for all $t\in[0,1]$,
\begin{equation}\label{e1}
\lambda_t\left\{ |\Phi_n(N_n)|-\beta \Phi(|N_n|)\right\} \le \frac{1}{4\beta}\lambda_t\{\Phi_n(|N_n|)\}.
\end{equation}

Since $\Phi_n(N_n)\to \Phi(N)$, we also have $|\Phi_n(N_n)|\to |\Phi(N)|$ thanks to Tykhonov's theorem \cite{T} (see also \cite[Theorem 1.1]{DDPS}), and $|N_n|\to |N|$. Therefore
\begin{equation}\label{e2}
   |\Phi_n(N_n)|-\beta \Phi(|N_n|) \to |\Phi(N)|-\beta \Phi(|N|) \quad {\mathrm{and}} \quad \Phi_n(|N_n|)\to \Phi(|N|)
\end{equation}
Combining \eqref{e1}, \eqref{e2}, and Lemma \ref{lemma-ae-convergence}, we obtain that for almost every $t\in(0,1)$,
\begin{equation*}
\lambda_t\left\{ |\Phi(N)|-\beta \Phi(|N|)\right\} \le \frac{1}{4\beta}\lambda_t\{\Phi(|N|)\}.
\end{equation*}
Since the function $t\mapsto \lambda_t(A)$ is right-continuous, this relation actually holds for all $t\in(0,1)$ and Lemma \ref{lemma-sp-dominance}
completes the proof.
\end{proof}

Recall that $\mathcal{M}$ denotes a semifinite von Neumann algebra and let ${\mathcal{F}}_{\infty}$ be a type ${\mathrm{II}}_{\infty}$ factor.

\vskip 5pt
\begin{theorem}\label{thminfinite}     Let $\Phi: L_p({\mathcal{M}})\to \overline{\mathcal{F}}_{\infty}$ be a continuous positive linear map, let $N\in L_p({\mathcal{M}})$ be normal, $\beta> 0$ and $\varepsilon>0$. Then, there exists a partial isometry $V\in{\mathcal{F}}_{\infty}$ such that
$$
   |\Phi(N)|  \le \beta\Phi(|N|) +\frac{1}{4\beta}V\Phi(|N|)V^* +\varepsilon I.
$$
\end{theorem}

\vskip 5pt
The proof is similar to the type ${\mathrm{II}}_1$ case, except that we use the following lemma instead of Lemma \ref{lemma-sp-dominance}  (Lemma \ref{lemma-ae-convergence} still holds in the type ${\mathrm{II}}_{\infty}$ case with $t\in(0,\infty)$).

\vskip 5pt
\begin{lemma}\label{lemma-sp-dominance2}
For two self-adjoint operators  $A,B\in\overline{\mathcal{F}}_{\infty}$, the spectral dominance relation
$\lambda_t(A)\le\lambda_t(B)$ for all $t\in(0,\infty)$ holds if and only if for every $\varepsilon>0$ there
exists a partial isometry $U\in{\mathcal{F}}_{\infty}$ such that $A\le UBU^*+\varepsilon I$.
\end{lemma}

Recall that the Russo-Dye theorem says that positive linear maps on a unital $C^*$-algebras attain their norms at the identity. A generalization of this theorem was given in Corollary \ref{corRD1}. This can also be generalized to possibly unbounded positive linear maps as in our next two corollaries. We state the corollaries for type  ${\mathrm{II}}_1$ factors and for a map defined on $\overline{\mathcal{M}}$, of course similar statements hold for type ${\mathrm{II}}_{\infty}$ and/or for maps defined on  $L_p({\mathcal{M}})$, $p>0$.

\vskip 5pt
\begin{cor}\label{corRD2}  Let $Z\in \mathcal{M}$ be a contraction, let $\Phi: \overline{\mathcal{M}}\to \overline{\mathcal{F}}_1$ be a positive linear map continuous with respect to the measure topologies, and let $\varepsilon>0$. Then, for some unitary $V\in{\mathcal{F}}_1$,
$$
   |\Phi(Z)| \le \frac{\Phi(I) +V\Phi(I)V^*}{2} +\varepsilon I.
$$
\end{cor}

\vskip 5pt
\begin{proof} We repeat the proof of Corollary \ref{corRD1}: Dilate $Z$ into a unitary $U\in\bM_2( \overline{\mathcal{M}})$,
$$
U=\begin{pmatrix}
Z&-\sqrt{I-ZZ^*} \\
\sqrt{I-Z^*Z}& Z^*
 \end{pmatrix},
$$
and define $\Psi:\bM_2( \overline{\mathcal{M}})\to\overline{\mathcal{F}}_1$ by
$$
\Psi\left( \begin{pmatrix} A&B \\ C&D\end{pmatrix}\right)=\Phi(A).
$$
Applying Theorem \ref{thmfinite} to $\Psi$ and $U$, we have
$$
|\Phi(Z)|=|\Psi(U)|\le \frac{\Psi(|U|) +V\Psi(|U|)V^*}{2}= \frac{\Phi(I) +V\Phi(I)V^*}{2}
$$
for some unitary $V\in{\mathcal{F}}_1$.
\end{proof}

\vskip 5pt
An immediate consequence of Corolloray \ref{corRD2} is the following generalized $s$-number estimate.

\vskip 5pt
\begin{cor}\label{corcontraction2}  Let $Z\in \mathcal{M}$ be a contraction and let $\Phi: \overline{\mathcal{M}}\to \overline{\mathcal{F}}_1$ be a positive linear map continuous with respect to the measure topologies. Then, for all $t\in(0,1/2)$,
$$
  \mu_{2t} (\Phi(Z)) \le \mu_t(\Phi(I)).
$$
\end{cor}

\vskip 5pt
We close the paper by the following version of Corollary \ref{corcons2} whose proof is quite similar.

\vskip 5pt
\begin{cor}\label{corfinal}     Let $X\in \overline{\mathcal{M}}$, let $\Phi: \overline{\mathcal{M}}\to \overline{\mathcal{F}}_1$  be a positive linear map continuous with respect to the measure topologies, and let $\varepsilon>0$. Then, for some unitary $V\in{\mathcal{F}}_1$,
$$
   |\Phi(X+X^*)| \le \Phi(|X|+|X^*|) + \frac{1}{4}V\Phi(|X|+|X^*|)V^* +\varepsilon I.
$$
\end{cor}

\vskip 5pt
\noindent
Jean-Christophe Bourin\\
Laboratoire de math\'ematiques,
Universit\'e de Franche-Comt\'e,
25 000 Besan\c{c}on, France.\\
Email: jcbourin@univ-fcomte.fr

\vskip 5pt
\noindent
Jingjing Shao\\
School of Mathematics and Statistic Sciences, Ludong University, Yantai 264001, China.\\
Email: jingjing.shao86@yahoo.com

\end{document}